\documentclass[11pt]{amsart}
\usepackage{amssymb}
\usepackage{layout} 

\newtheorem{theorem}{Theorem}

\newtheorem{corollary}{Corollary}

\newcommand{\Z}{\mathbb Z}

\newcommand{\U}{\mathbb U}

\begin{document}

\title[Exponential Sums]{exponential Sums and Distinct Points on Arcs} 
\author[\O. J. R\o dseth]{\O ystein J. R\o dseth}
\address{Department of Mathematics, University of Bergen, Johs.
  Brunsgt. 12, N-5008 Bergen, Norway} 
\email{rodseth@math.uib.no}
\urladdr{http://math.uib.no/People/nmaoy}

\begin{abstract}
Suppose that some harmonic analysis arguments have been invoked to show that the indicator function of a set of residue classes modulo some integer has a large Fourier coefficient. To get information about the structure of the set of residue classes, we then need a certain type of complementary result. A solution to this problem was given by Gregory Freiman in 1961, when he proved a lemma which relates the value of an exponential sum with the distribution of summands in semi-circles of the unit circle in the complex plane. Since then, Freiman's Lemma has been extended by several authors. Rather than residue classes, one has considered the situation for finitely many arbitrary points on the unit circle. So far, Lev is the only author who has taken into consideration that the summands may be bounded away from each other, as is the case with residue classes. In this paper we extend Lev's result by lifting a recent result of ours to the case of the points being bounded away from each other.
\end{abstract}

\dedicatory{Dedicated in honour of Mel Nathanson's $60^{th}$ birthday\\{\rm (Contribution to the Festschrift)}}
\date{\today}
\maketitle

\section{Introduction} \label{sec:intro}
In additive combinatorics, and in additive combinatorial number theory in particular, situations of the following type are rather common. Let $A$ be a set of $N$ residue classes modulo an integer $m$. Suppose that some harmonic analysis arguments have been invoked to show that the indicator function of $A$ has a large Fourier coefficient; i.e.,
\[
 \max_{0\neq x\in\Z/m\Z}|\widehat{1}_A(x)|\geq\alpha N
\]
for some $\alpha\in(0,1]$. Following Green \cite{gre}, or see Tao and Vu \cite{t&v}, we might say that $A$ has ``Fourier bias''. Using this information one wishes to conclude that $A$ has ``combinatorial bias''; perhaps integer representatives of the residue classes in some affine image of $A$ concentrates on some interval. For this we then need  tight upper bounds for the absolute values of the exponential sums $\widehat{1}_A(x)$. Now, the basic idea is that the absolute value of an exponential sum is small if the terms are, in some sense, uniformly distributed. The case $m$ prime was applied by Freiman \cite{fre3,fre4} in the proof of his ``2.4-theorem''. A proof of this theorem is presented in both Freiman's  classical monograph \cite{fre} and in Mel Nathanson's beautiful book \cite{mel}.

Freiman \cite{fre3,fre4}, Postnikova \cite{pos}, Moran and Pollington \cite{m&p}, Lev \cite{lev2,lev3}, and R\o dseth \cite{rod1,rod2}  considered the situation in which one has finitely many arbitrary points on the unit circle in the complex plane, rather than a subset of $\Z/m\Z$. In possible applications the points on the unit circle will sometimes be bounded away from each other, like, for instance, if we are looking at a subset of $\Z/m\Z$. If we add the condition that the points should be bounded away from each other, we could hope for sharper results. Indeed, by adding this assumption, Lev sharpened Freiman's Lemma to \eqref{bound1} below. Lev's result \cite[Theorem~2]{lev2} seems, however, to be the only result in the literature addressing this issue. In this paper we prove a result which extends both Lev's result about points  being bounded away from each other, and our main result in \cite{rod1,rod2}, where we did not take this property into consideration.

First, we shall, however, present a version of Freiman's Lemma. We include two different proofs in the hope of giving the reader an impression of two recent techniques used in the search and study of results related to Freiman's Lemma.

\section{Three Theorems}  \label{sec_results}
\noindent
In the following, $n$, $N$, $\kappa$, and $k$ are non-negative or positive integers. We write $\U$ for the unit circle in the complex plane; that is, 
\[
\U:=\{z\in\mathbb C :|z|=1\}.
\]
An empty sum is taken as zero.

We now state Freiman's Lemma \cite{fre3,fre4}.

\begin{theorem}[Freiman's Lemma] \label{FL}
Suppose that the complex numbers $z_1,\dots,z_N\in\U$ have the property that any open semi-circle of\, $\U$ contains at most $n$ of them. Then
\begin{equation} \label{fl} 
|z_1+\cdots+z_N|\leq2n-N.
\end{equation}
\end{theorem}

The complex numbers $z_1,\dots,z_N$ are not necessarily all distinct.
The assumptions of the lemma imply that $N\leq2n$, and the result is
{\em sharp} in the range $n\leq N\leq2n$. That is, for every $n$ and
$N$ satisfying these inequalities, there exist sequences $z_1,\ldots,z_N$
with $z_j\in\U$, such that the hypotheses are satisfied and
$|z_1+\dots+z_N|$ meets the bound $2n-N$.  In this sense, Freiman's
Lemma is best possible. The result has, however, been extended by Moran and Pollington \cite{m&p}, Lev \cite{lev2, lev3}, and by R\o dseth \cite{rod1,rod2}.

The next theorem is Lev's sharpening of Freiman's Lemma at the expense of requiring the points to be bounded away from each other. Lev \cite[Theorem~2]{lev2} proved the following theorem.

\begin{theorem}[Lev] \label{thm1}
Let $\delta\in(0,\pi]$ satisfy $n\delta\leq\pi$. Suppose that the complex numbers $z_1,\dots,z_N\in\U$ have the following two properties:
\begin{enumerate}
 \item[(a)] any open semi-circle of $\,\U$ contains at most $n$ of them;
\item[(ii)] any open arc of $\,\U$ of length $\delta$ contains at most one of them.
\end{enumerate}
Then
\begin{equation} \label{bound1}
 |z_1+\dots+z_N|\leq\frac{\sin((2n-N)\delta/2)}{\sin(\delta/2)}.
\end{equation}
\end{theorem}

So, by introducing the condition (ii), Lev reduced Freiman's bound $2n-N$ in \eqref{fl} to the bound in \eqref{bound1}, a refinement that, according to Lev, was crucial in \cite{lev4}. The result is sharp for $n\leq N\leq2n$. By letting $\delta\to0^+$ in Lev's result, we recover Freiman's Lemma.

In this paper we shall prove the following theorem.
\begin{theorem} \label{main}
Let $\delta,\varphi\in(0,\pi]$ satisfy $n\delta\leq\varphi$. Suppose that the complex numbers $z_1,\dots,z_N\in\U$ have the following two properties:
\begin{enumerate}
 \item[(i)] any open arc of $\,\U$ of length $\varphi$ contains at most $n$ of them;
\item[(ii)] any open arc of $\,\U$ of length $\delta$ contains at most one of them.
\end{enumerate}
Let $N=\kappa n+r$, $1\leq r\leq n$, and assume that $(\kappa+1)\varphi\leq2\pi$. Then we have
\begin{equation} \label{th2}
|z_1+\cdots+z_N|
\leq\frac{\sin(r\delta/2)}{\sin(\delta/2)}\cdot
\frac{\sin((\kappa+1)\varphi/2)}{\sin(\varphi/2)}
+\frac{\sin((n-r)\delta/2)}{\sin(\delta/2)}\cdot
\frac{\sin(\kappa\varphi/2)}{\sin(\varphi/2)}.
\end{equation}
\end{theorem}

Notice that $\kappa=\lceil N/n\rceil-1$. Also notice that we do need some condition like $(\kappa+1)\varphi\leq2\pi$, to be certain that there is room enough on $\U$ to place $N$ points such that they satisfy the other conditions set in the theorem. This condition can also be written as $N\leq\lfloor2\pi/\varphi\rfloor n$.

The restriction $n\delta\leq\varphi$ is no problem. For if $n\delta\geq\varphi$, then (i) follows from (ii) and can be omitted. Now, $|z_1+\cdots+z_N|$ attains its maximum on any $N$-term geometric progression with ratio $\exp(i\delta)$; hence
\begin{equation} \label{extra}
|z_1+\dots+z_N|\leq\left|\sum_{j=0}^{N-1} \exp(ij\delta)\right|=\frac{\sin(N\delta/2)}{\sin(\delta/2)}.
\end{equation}

The upper bound \eqref{th2} is attained on the union of two finite $2$-dimensional geometric progressions on $\U$, one consisting of geometric progressions with ratio $\exp(i\delta)$ and $r$ terms each, centered around the points
\[
\exp((-\kappa+2j)i\varphi/2),\quad j=0,1,\dots,\kappa,
\]
and the other consisting of geometric progressions with ratio $\exp(i\delta)$ and $n-r$ terms each, centered around the points
\[
\exp((-\kappa+1+2j)i\varphi/2),\quad j=0,1,\dots,\kappa-1.
\]
This shows that Theorem~\ref{main} is sharp for $N\leq\lfloor2\pi/\varphi\rfloor n$.

Clearly, the bound \eqref{th2} in Theorem~\ref{main} can be replaced by the weaker, but smooth and nice, bound
\[
|z_1+\cdots+z_N|\leq\frac{\sin(n\delta/2)}{\sin(\delta/2)}\cdot\frac{\sin(\varphi N/(2n))} {\sin(\varphi/2)};
\]
cf. Lev \cite{lev3}.

\section{Two Proofs of Freiman's Lemma} \label{sec:two}
For real numbers $\alpha<\beta$, the set of $z\in\U$ satisfying $\alpha<\arg z\leq\beta$ for some value of $\arg z$, is an \textit{open-closed} arc. In Freiman's Lemma one often assumes that any open-closed (or closed-open) semi-circle of $\U$ contains at most $n$ of the points $z_j$, instead of taking open semi-circles; cf. \cite{fre}, \cite{mel}. It is, however, easy to see that the two variants of the hypotheses are equivalent; cf. Section~\ref{arcs} below.
 
Freiman's proof of Theorem~\ref{FL} was simplified by Postnikova \cite{pos}, and it is this proof we find in the books \cite{fre} and \cite{mel}. Here, we shall present two other proofs in an attempt to give the reader a pleasant introduction to two  techniques recently employed in the quest for extensions of Freiman's Lemma. The two proofs are rather different. One proof can be characterized as topological-combinatorial or as a perturbation method, and is due to Lev (extracted from the proof of \cite[Theorem~2]{lev2}). The other proof uses properties of a certain Fourier coefficient, and is independently due to Lev and the present author.

\subsection{First Proof}
Assume that Freiman's Lemma is false. We consider the smallest $N$ for which there exists an $N$-term sequence which satisfies the hypotheses, but violates \eqref{fl}. Then $N>1$. By considering open semi-circles, the set of $N$-term sequences satisfying the hypotheses form a closed subset of the compact topological space $\U^N$, and is itself compact. By the continuity of the function $z_1,\ldots,z_N\mapsto|S|$, where $S:=z_1+\cdots+z_N$, we thus have that $|S|$ attains a maximum value on some  $N$-term sequence $Z:=z_1,\ldots,z_N$, which satisfies the hypotheses. Then $|S|>2n-N$. A rotation of $Z$ shows that we may assume that $S$ is real and non-negative.

Suppose that there is a $z\neq1$ in $Z$. By symmetry, we may assume that $\mathrm{Arg}\,z<0$, using the interval $[-\pi,\pi)$ for the principal argument.  Replacing $z$ by $z\exp(i\varepsilon)$ for a small $\varepsilon>0$, we get an increase in $|S|$. Thus the replacement results in violation of the hypotheses in Freiman's Lemma. The only possibility is that the replacement produces an open semi-circle with more than $n$ points from $Z$; hence $-z$ also belongs to $Z$.

We now remove $\pm z$ from $Z$. This gives us a sequence $Z'$ satisfying the hypotheses, and with parameters $N'=N-2$ and $n'=n-1$. Denoting the sum of the terms of $Z'$ by $S'$, we have by the minimality of $N$,
\[
 S=S'\leq2n'-N'=2n-N,
\]
a contradiction. Thus all terms of $Z$ are equal to $1$, and $N=S>2n-N$. But a semi-circle containing $1$ contains $N$ terms from $Z$; hence $n\geq N$, and again we have a contradiction.

\subsection{Second Proof}
Assume that $z_1,\ldots,z_N$ satisfy the hypotheses. We shall prove \eqref{fl}, and can without loss of generality assume that $S$ is real and non-negative. Let $K(\theta)$ denote the number of values $j\in[1,N]$ such that $\theta-\pi\leq\arg z_j<\theta$ for some value of $\arg z_j$. Then we have
\begin{equation} \label{K}
K(\theta)+K(\theta+\pi)=N,
\end{equation}
so that $N-n\leq K(\theta)\leq n$.

Moreover,
\[
\int_{-\pi}^\pi K(\theta)\sin\theta\,d\theta=\sum_{j=1}^N\int_{\arg z_j}^{\arg z_j+\pi}\sin\theta\,d\theta=2S,
\]
and we obtain
\begin{eqnarray*}
2S&=&\left(\int_{-\pi}^0+\int_0^{\pi}\right)K(\theta)\sin\theta\,d\theta\\
{}&\leq& (N-n)\int_{-\pi}^0\sin\theta\,d\theta+n\int_0^{\pi}\sin\theta\,d\theta\\
{}&=&4n-2N.
\end{eqnarray*}

\section{Proof of Theorem~\ref{main}} \label{sec:proof2}
We now turn to the proof of Theorem~\ref{main}. In an attempt to make the proof more readable, we split the proof up into several parts.

\subsection{Notation}
Set 
\[
\omega=\exp(i\varphi/2)\quad \text{and}\quad\rho=\exp(i\delta/2).
\]
We often denote a sequence $z_1,\ldots,z_N\in\U^N$ by $Z$, and we write $S=z_1+\cdots+z_N$. If the sequence $Z\in\U^N$ satisfies the assumptions of Theorem~\ref{main} for a certain value of $n$, we say that $Z$ is an $(N,n)$-\textit{admissible} sequence. We use the interval $[-\pi,\pi)$ for the principal argument $\mathrm{Arg}\,z$ of a non-zero complex number $z$. 

\subsection{Arcs}  \label{arcs}
An arc $(u,v)$, where $u,v\in\U$, consists of the points we pass in moving counter-clockwise from $u$ to $v$. Almost all arcs in this paper have lengths at most $2\pi$.
In spite of the similarity of notation, an arc cannot be confused with a real interval.

Consider the two statements
\begin{enumerate}
 \item[(i)] any open arc of $\,\U$ of length $\varphi$ contains at most $n$ of the points $z_j$;
\item[(i$'$)] any open-closed arc of $\,\U$ of length $\varphi$ contains at most $n$ of the points $z_j$.
\end{enumerate}
Clearly, (i$'$) implies (i). On the other hand, if (i$'$) fails, then (i) fails. For if an open-closed arc $(u,u\exp(i\varphi)]$, $u\in\U$, of $\U$ contains $n+1$ of the points $z_j$, then, for a sufficiently small real $\varepsilon>0$, the open arc $(u\exp(i\varepsilon),u\exp(i\varphi+i\varepsilon))$ contains the same $n+1$ points. Thus (i) and (i$'$) are equivalent. (And the reason is, of course, that we only consider finitely many points $z_j$.)

\subsection{Assumptions}
The two bounds \eqref{th2} and \eqref{extra} coincide for $n\delta=\varphi$. For the proof of Theorem~\ref{main} we may therefore assume that $n\delta<\varphi$. Moreover, we observe that the right-hand side of \eqref{th2} is continuous as a function of $\varphi$ on the real interval $(n\delta,2\pi/(\kappa+1)]$. Hence it suffices to prove the assertion of Theorem~\ref{main} in the case $\varphi<2\pi/(\kappa+1)$. 

We shall prove Theorem~\ref{main} by contradiction. We therefore assume the theorem false. Choose the smallest non-negative integer $N$ for which there exists an $n$ such that \eqref{th2} fails for some $(N,n)$-admissible sequence. Then $N>1$.

\subsection{Geometric Progressions}
A geometric progression $\Delta$ in $\U$ with ratio $\rho^2$ is called a $\delta$-\textit{progression}. If we write this progression as 
\begin{equation} \label{del}
u\rho^{r-1-2j},\quad j=0,\ldots,r-1,
\end{equation}
for some $u\in\U$, then $\Delta$ is a progression of \textit{length} $r$, \textit{centered around} $u$. The point $u$ may, or may not, belong to the progression. The point $u\rho^{r-1}$ is the \textit{first} and $u\rho^{-(r-1)}$ is the \textit{end} point (element, term) of the progression. If we multiply each term of $\Delta$ by $v\in\U$, we get a new $\delta$-progression denoted $\Delta v$. If all terms of $\Delta$ belong to $Z$, we say that $\Delta$ is \textit{in} $Z$. The $\delta$-progression \eqref{del} is \textit{maximal} in $Z$, if the progression is in $Z$, but neither $u\rho^{r+1}$ nor $u\rho^{-r-1}$ belongs to $Z$.

Notice that if $\Delta$ is a $\delta$-progression in $Z$ of length $r$, then $r\leq n$. The reason is the inequality $n\delta<\varphi$, and that an open arc of $\U$ of length $\varphi$ contains at most $n$ points from $Z$.

\subsection{Compactness}
The arcs in both (i) and (ii) are open. Therefore, by a standard compactness argument, $|S|$ attains a maximum on  the set of $(N,n)$-admissible sequences. Let the maximum be attained at the $(N,n)$-admissible sequence $Z:=z_1,\ldots,z_N$. A rotation of $Z$ shows that we may assume that the matching $S=z_1+\cdots+z_N$ is real and non-negative.

\subsection{Dispersion}
Recall that $N>1$ and $N=\kappa n+r$, $1\leq r\leq n$. Suppose that $Z$ is a $\delta$-progression. The length of this $\delta$-progression is $N$.

If $\kappa=0$, then
\[
 |z_1+\cdots+z_N|\leq\frac{\sin(r\delta/2)}{\sin(\delta/2)},
\]
and \eqref{th2} holds, contrary to hypothesis. Thus $N\geq n+1$, and there is a closed arc of $\U$ of length $n\delta$ containing $n+1$ of the points $z_j$. But we have $n\delta<\varphi$, so this contradicts condition (i) of Theorem~\ref{main}. Therefore $Z$ is not a $\delta$-progression.

\subsection{Perturbation} \label{sec:per}
We have just seen that $Z$ is not a $\delta$-progression. This implies that there exists at least one point $z\in Z$, $z\neq0$, such that $\mathrm{Arg}\,z>0$ and $z\rho^{-2}\not\in Z$, or $\mathrm{Arg}\,z<0$ and $z\rho^2\not\in Z$. Choose such a $z$, as close to  $-1$ as possible. By symmetry, we can assume that $\mathrm{Arg}\,z<0$. We split this case into the two cases determined by \eqref{ulik} and \eqref{venstre} below.

\subsubsection{Case I}
First, we assume that there is an integer $\lambda$ in the interval $1\leq\lambda\leq\kappa$ satisfying
\begin{equation} \label{ulik}
-\frac{\lambda\varphi}{2}\leq\mathrm{Arg}\,z<-\frac{(\lambda-1)\varphi}{2}.
\end{equation}

The length of the shortest arc $(z,z')$ for $z\neq z'\in Z$ is at least $\delta$. Since $z\rho^2\not\in Z$, the arc $(z,z')$ has length greater than $\delta$. 
We start the perturbation procedure by performing a small counter-clockwise rotation of $z$ along the unit circle; that is, we replace $z$ by  $z\exp(i\varepsilon)$ for a small $\varepsilon>0$. If $\varepsilon$ is sufficiently small, then the length of the open arc $(z\exp(i\varepsilon),z')$ is still at least $\delta$, and the rotation of $z$ does not disturb the truth of (ii). Since $\mathrm{Arg}\,z<0$, there will, however, be an  increase in $|S|$. Thus the rotation violates (i$'$). Therefore we have $z\omega^2\in Z$, and the open-closed arc $(z,z\omega^2]$ contains exactly $n$ terms from $Z$. If $z\omega^2\rho^2\in Z$ then the arc $(z\rho^2,z\omega^2\rho^2]$ contains $n+1$ points of $Z$. Therefore $z\omega^2\rho^2\not\in Z$.

Next, we perform a small rotation of both the points $z$ and $z\omega^2$ simultaneously and  counter-clockwise along the unit circle. We have
\[
\mathrm{Arg}(z+z\omega^2)=\mathrm{Arg}\,z+\frac{\varphi}{2},
\]
and if $\lambda\geq2$, the rotation makes $|S|$ increase. Hence $z\omega^4\in Z$. We also see that $z\omega^4\rho^2\not\in Z$, and that the arc $(z\omega^2,z\omega^4]$ contains exactly $n$ terms from $Z$. 

Using the right inequality in \eqref{ulik}, we find for $1\leq j\leq \lambda$, that
\[
\mathrm{Arg}\left(\sum_{\ell=0}^{j-1} z\omega^{2\ell}\right)=\mathrm{Arg}\,z+\frac{(j-1)\varphi}{2}<-\frac{(\lambda-j)\varphi}{2}\leq0.
\]
This shows that we may continue the perturbation process until we eventually obtain
\[
z\omega^{2j}\in Z\quad\text{and}\quad z\omega^{2j}\rho^2\not\in Z\quad\text{for}\quad j=0,1,\ldots,\lambda.
\]
In addition, for $j=1,2,\ldots,\lambda$, each arc $(z\omega^{2(j-1)},z\omega^{2j}]$ contains exactly $n$ points from $Z$.

By the left inequality in \eqref{ulik}, we have
\[
\mathrm{Arg}(z\omega^{2\lambda})\geq-\mathrm{Arg}\,z.
\]
Using the definition of $z$, it follows that if we start at the point $z\omega^{2\lambda}\rho^2$, which is in $\U$ but not in $Z$, and move counter-clockwise on $\U$, then the first point in $Z$ we meet, is the end point of a maximal $\delta$-progression in $Z$ with $z$ as its first element.

The open-closed arc $(z,z\omega^{2\lambda}]$ contains $\lambda n$ points from $Z$, and in addition we have $r'\geq1$ points in the maximal $\delta$-progression in $Z$ with $z$ as first point. Thus we have
\[
\lambda n+r'=N=\kappa n+r, \quad 1\leq r',r\leq n,
\]
so that $\lambda=\kappa$ and $r'=r$. This completes the first case \eqref{ulik}.

\subsubsection{Case II}
Second, if $z$ does not satisfy \eqref{ulik} for $\lambda=\kappa$, then
\begin{equation}  \label{venstre}
-\pi\leq\mathrm{Arg}\,z<-\frac{\kappa\varphi}{2},
\end{equation}
and we can do one more step in the perturbation process. By the inequality $(\kappa+1)\varphi<2\pi$, we have that the arc $(z,z\omega^{2(\kappa+1)}]$ contains exactly $(\kappa+1)n$ points from $Z$. In addition, we have the point $z$ itself. Thus we have found $(\kappa+1)n+1$ points in $Z$. Since $Z$ contains at most $(\kappa+1)n$ points, we have a contradiction. Thus there are no points in $Z$ satisfying \eqref{venstre}.

\subsection{Primary Points}
For the $z\in Z$ defined at the beginning of Section~\ref{sec:per}, we now know that for $j=1,2,\ldots,\kappa$, each open-closed arc $(z\omega^{2(j-1)},z\omega^{2j}]$ contains exactly $n$ points from $Z$. In addition, there is the maximal $\delta$-progression $\Delta$ of length $r$ with $z$ as first element. This accounts for all elements in $Z$.

Thus we have
\[
z\omega^{2j}\in Z,\quad j=0,1,\ldots,\kappa.
 \]
Suppose that $z\rho^{-2}\in\Delta$. The open-closed arc $(z,z\omega^2]$ contains exactly $n$ points from $Z$. If $z\omega^2\rho^{-2}\not\in Z$, then the closed-open arc $[z\rho^{-2},z\omega^2\rho^{-2})$ contains $n+1$ points from $Z$, contrary to hypotheses. Thus $z\omega^2\rho^{-2}\in Z$. Continuing in the same manner, we get (with a slight abuse of notation) $\Delta\omega^2\subseteq Z$. Repeating the argument, we ultimately obtain,
\[
 \Delta\omega^{2j}\subseteq Z\quad\text{for $j=0,1,\ldots,\kappa$}.
 \]
Notice that for $j\geq1$, the $\delta$-progression $\Delta\omega^{2j}$ is not necessarily maximal.  

We set
\[
Z_1=\bigcup_{j=0}^\kappa \Delta\omega^{2j}.
\]
Then $N_1:=|Z_1|=\kappa r+r$. We write $S_1$ for the sum of all elements in $Z_1$, and we want to find an upper bound for $|S_1|$. (We may have $N_1=N$, so we cannot use \eqref{th2}.) We can, however, easily determine the exact value of $|S_1|$,
\begin{equation} \label{s1}
|S_1|=\left|\sum_{\ell=0}^{r-1}\rho^{-2\ell}\sum_{j=0}^\kappa z\omega^{2j}\right|
=\frac{\sin(r\delta/2)}{\sin(\delta/2)}
\cdot\frac{\sin((\kappa+1)\varphi/2)}{\sin(\varphi/2)}.
\end{equation}
If $\kappa=0$ or $r=n$, then $Z=Z_1$, and by \eqref{s1}, Eq. \eqref{th2} holds, contrary to hypothesis. Thus we have $\kappa\geq1$ and $1\leq r<n$.

\subsection{Secondary Points}
We set
\[
Z_2=Z\setminus Z_1,
\]
and put $N_2=|Z_2|$. Then $N_2=(\kappa-1)(n-r)+(n-r).$ We have that $Z_2$ is $(N_2,n-r)$-admissible, and since $N_2<N$, we can use \eqref{th2}. Writing $S_2$ for the sum of the terms in $Z_2$, we obtain
\begin{equation} \label{s2}
|S_2|\leq\frac{\sin((n-r)\delta/2)}{\sin(\delta/2)}
\cdot\frac{\sin(\kappa\varphi/2)}{\sin(\varphi/2)}.
\end{equation}
We have
\[
S\leq|S_1|+|S_2|,
\]
and \eqref{s1} and \eqref{s2} show that \eqref{th2} holds; again we have reached a contradiction. This concludes the proof of Theorem~\ref{main}.

\section{Closing Remarks}
\label{sec:turn}
We can also state Theorem~\ref{main} in a slightly different way.

\begin{theorem} \label{altmain}
Let $\delta,\varphi\in(0,\pi]$ satisfy $n\delta\leq\varphi$. Assume that the complex numbers $z_1,\dots,z_N\in\U$ have the two properties {\normalfont(i)} and {\normalfont(ii)} stated in Theorem~\ref{main}. Then we have
\begin{multline} \label{alt}
|z_1+\cdots+z_N|
\leq\frac{\sin((N-(k-1)n)\delta/2)}{\sin(\delta/2)}\cdot
\frac{\sin(k\varphi/2)}{\sin(\varphi/2)}\\
+\frac{\sin((kn-N)\delta/2)}{\sin(\delta/2)}\cdot
\frac{\sin((k-1)\varphi/2)}{\sin(\varphi/2)},
\end{multline}
for any positive integer $k\leq2\pi/\varphi$.
\end{theorem}

Denote the right-hand side of \eqref{alt} by $L_k$. We easily see that
\[
L_{k+1}-L_k=2\frac{\sin((kn-N)\delta/2)}{\sin(\delta/2)}\cdot\frac{\sin(k\varphi/2)}{\sin(\varphi/2)}\left(\cos\frac{n\delta}{2}-\cos\frac{\varphi}{2}\right),
\]
so that
\[
L_1\geq\ldots\geq L_{\lceil N/n\rceil}\leq L_{\lceil N/n\rceil+1}\leq\ldots\leq L_{\lfloor2\pi/\varphi\rfloor}.
\]
Thus $L_k$ attains its minimum at $k=\lceil N/n\rceil$ (and also at $k=N/n+1$ if $n$ divides $N$).  Therefore Theorem~\ref{main} and Theorem~\ref{altmain} are essentially one and the same.

Now, by letting $\delta\to0^{+}$, we recover the main result of \cite{rod1,rod2}.

The following direct extension of Theorem~\ref{thm1} is an immediate consequence of  Theorem~\ref{altmain}.

\begin{corollary} \label{cor}
Let $k\geq2$ be an integer, and let $\delta\in(0,\pi]$ satisfy $n\delta\leq2\pi/k$. Suppose that the complex numbers $z_1,\dots,z_N\in\U$ have the following two properties:
\begin{enumerate}
 \item[(b)] any open arc of $\,\U$ of length $2\pi/k$ contains at most $n$ of them;
\item[(ii)] any open arc of $\,\U$ of length $\delta$ contains at most one of them.
\end{enumerate}
Then
\[ 
 |z_1+\dots+z_N|\leq\frac{\sin((kn-N)\delta/2)}{\sin(\delta/2)}.
\]
\end{corollary}

If $k$ is even, then Corollary~\ref{cor} is an easy consequence of Theorem~\ref{thm1}. This does not seem to be the case for $k$ odd. Corollary~\ref{cor} is sharp for $(k-1)n\leq N\leq kn$. By letting $\delta\to0^+$, we recover the result of Moran and Pollington \cite{m&p}.

For applications, one would perhaps, although not always necessary, turn the upper bound for $|z_1+\cdots+z_N|$ into a lower bound for the number of terms $z_j$ in at least one open arc of $\U$ of length $\varphi$.

In closing, let us consider such a lower bound for the residue class situation mentioned in the introduction. Then we have an $N$-set $A$ of residue classes modulo an integer $m>1$.  Set $\delta=2\pi/m$, and put $\varphi=2\pi/k$ for some integer $k\geq2$.

Let
\[
\widehat{1}_A(1)=\sum_{a\in A}\exp(2\pi ia/m).
\]
By Corollary~\ref{cor}, there exist integers $u,v$ satisfying $u\leq v<u+m/k$ such that the image of the interval $[u,v]$ under the canonical homomorphism $\Z\to\Z/m\Z$ contains $n_0$ elements of $A$, where
\[
n_0\geq\left\lceil\frac{1}{k}\left(N+\frac{\mathrm{Arcsin}(|\widehat{1}_A(1)|   \sin(\pi/m))}{\pi/m}\right) \right\rceil;
\]
cf. \cite[Corollary~2]{lev2}, \cite[Section~6]{rod1}.

\end{document}